\PassOptionsToPackage{hidelinks}{hyperref}
\documentclass{siamart220329}
\usepackage{graphicx}
\usepackage{xcolor}
\usepackage{textcomp}
\usepackage{booktabs}
\usepackage{multirow}
\usepackage{mathrsfs}
\usepackage{algorithm}
\usepackage{algorithmic}   
\usepackage{listings}
\usepackage{subfig}        
\usepackage{float}         
\usepackage{placeins}
\usepackage{adjustbox}
\usepackage{siunitx}
\usepackage{pdfpages}
\usepackage{pdflscape}
\usepackage{tabularx}
\usepackage{titlesec}

\titlespacing*{\section}{0pt}{1ex}{0.8ex} 

\title{Efficient Multi-Precision Computation of Bessel Functions for Real Orders and Complex Arguments with Fortran Implementation--- Part~II: The Modified Bessel Function of the Second Kind, $K_\nu(z)$}

\author{
  Mofreh R.~Zaghloul%
  \thanks{Department of Physics, College of Sciences,
    United Arab Emirates University, Al~Ain 15551, Abu Dhabi, UAE
    (\email{m.zaghloul@uaeu.ac.ae}).}
  \thanks{Department of Mathematics,
    Massachusetts Institute of Technology, Cambridge, MA, USA.}
  \and
  Steven G.~Johnson%
  \thanks{Department of Mathematics,
    Massachusetts Institute of Technology,
    182 Memorial Dr., Cambridge, MA 02139, USA
    (\email{stevenj@math.mit.edu}).}
}

\headers{Efficient Computation of Modified Bessel Functions}{M. R.~Zaghloul and S. G.~Johnson}

\ifpdf
\hypersetup{
  pdftitle={Efficient and Accurate Computation of Modified Bessel Functions},
  pdfauthor={Mofreh R.~Zaghloul and Steven G.~Johnson}
}
\fi

\usepackage{amsmath,amssymb,amsfonts}
\begin{document}

\maketitle

\begin{abstract}
This paper, the second in a series, presents an efficient, self-contained algorithm for computing the modified Bessel function of the second kind, \(K_{\nu}(z)\) for complex argument and real orders, building on Part~I (for \(I_{\nu}\)). The method adaptively selects among analytic representations such as power series, large-\(|z|\) asymptotics, uniform asymptotics for large \(|\nu|\), and numerically stable forward recurrence with region boundaries tuned for accuracy and efficiency. A robust \texttt{Fortran} implementation supports double precision and quadruple precision. The use of quadruple precision extends the reliable computational domain and improves stability in challenging regimes. Accuracy is validated against high-precision \texttt{Maple} results, and benchmarks show runtimes significantly superior to those of established methods, in the literature, while avoiding their numerical failure modes across several decades of the parameter domain. Together with Part~I, this work provides a comprehensive, multiple-precision toolkit for \(\{I_{\nu},K_{\nu}\}\) across wide parameter ranges.
\end{abstract}



\section{Introduction}\label{sec:intro}
The modified Bessel functions of the first and second kinds, \( I_{\nu}(z) \) and \( K_{\nu}(z) \), play a fundamental role in numerous areas of applied mathematics, physics, and engineering. These functions naturally arise as solutions to the modified Bessel differential equation:
\begin{equation}
z^2 \frac{d^2 y}{dz^2} + z \frac{dy}{dz} - (z^2 - \nu^2) y = 0,
\end{equation}
where $\nu$ is the order, and $z \in \mathbb{C}$ is the complex argument. The order $\nu$ may, in general, be complex; however, in this study, we consider only the case of real order, $\nu \in \mathbb{R}$.

Accurate and efficient numerical evaluation of $I_\nu(z)$ and $K_\nu(z)$ is essential not only for direct use in physics and engineering applications, but also as a foundation for computing the broader family of Bessel and related functions. However, this task is computationally challenging due to the functions’ exponential behavior, branch points, and susceptibility to overflow and underflow in finite-precision arithmetic. The difficulties become more pronounced for complex arguments or large orders $\nu$, where different asymptotic and recurrence behaviors dominate and must be handled with care.

In the first part~\cite{Zaghloul_Johnson_2025}, we presented a robust and highly accurate algorithm for evaluating the modified Bessel function of the first kind, $I_\nu(z)$, over a wide domain in the complex plane with positive and negative real values of $\nu$. That work employed adaptive algorithmic switching between series expansion, asymptotic expansions, recurrence, and uniform asymptotic expansion to ensure both stability and efficiency across regimes. The present work is a natural continuation of that effort, extending our high performance and multiple-precision strategy to the computation of $K_\nu(z)$.

The modified Bessel function $K_\nu(z)$ has the integral representation:
\begin{equation}
   K_\nu(z) = \int_0^\infty e^{-z \cosh t} \cosh(\nu t)\, dt, \quad \text{Re}(z) > 0.
\end{equation}

For complex arguments \( z = x + iy \), this expression can be decomposed into two real-valued integrals:
\begin{equation}
K_{\nu}(z) = \mathcal{I}_{\cos} - i  \mathcal{I}_{\sin},
\end{equation}

where
\begin{equation}
\mathcal{I}_{\cos}(\nu, z) = \int_0^{\infty} e^{-x \cosh t} \cosh(\nu t) \cos(y \cosh t) \, dt, \quad \textit{x} > 0 . 
\end{equation}

\begin{equation}
\mathcal{I}_{\sin}(\nu, z) = \int_0^{\infty} e^{-x \cosh t} \cosh(\nu t) \sin(y \cosh t) \, dt,  \quad \textit{x} > 0.  
\end{equation}
These integrals may be evaluated using standard adaptive quadrature methods in double or extended precision. However, such approaches are generally slow and suffer from accuracy limitations in challenging regimes.

Several dedicated algorithms have been developed over the years to compute \( K_\nu(z) \) ~\cite{Temme1975, Campbell_81, Amos1983a, Amos1983b, Amos1986, Thompson_Barnett_1987, Press_Teukolsky_1991}. Most notably, Algorithm 644 by Amos ~\cite{Amos1986} is implemented in widely used libraries such as AMOS, SLATEC, and the GNU Scientific Library (GSL). While reliable within the confines of standard double precision arithmetic, Algorithm 644 enforces overly conservative overflow and underflow thresholds, leading to premature termination and loss of reliability in regions where valid results could still be computed. Moreover, it does not support higher precision, a feature increasingly needed in high-accuracy applications such as plasma physics, computational electromagnetics, astrophysical modeling, and quantum field theory. The need for and importance of calculating special functions in quadruple-precision arithmetic, as a compromise between the efficient but less accurate double-precision and the highly accurate but significantly slower variable-precision arithmetic offered by a few symbolic software packages, has been discussed in several previous publications~\cite{Zaghloul2023a,  Zaghloul2023b, Zaghloul2024, Zaghloul_Johnson_2025}.

It is worth noting that Algorithm~644, despite its historical importance, returns incorrect results when invoked for \(K_{\nu}(z)\) with negative orders (\(\nu<0\)). This is unexpected because \(K_{\nu}(z)\) is even in the order, satisfying \(K_{-\nu}(z)=K_{\nu}(z)\). In practice, this limitation can be easily circumvented by passing the absolute value of the order whenever \( \nu < 0\) , which enforces the identity and restores correctness.\\

In this work, we develop a new, numerically stable, and performance-efficient algorithm for computing \( K_\nu(z) \), implemented in modern, modular \texttt{Fortran} code supporting both double and quadruple (quad) precision. For small arguments, we employ power series and upward recurrence relations; for large arguments and/or large orders, we make use of asymptotic expansions. In the transition regions and for complex \( z \), we employ Temme’s well-received method, based on Miller's algorithm for evaluating a related function, the same foundational technique used in Algorithm 644, but integrated within a more flexible, more efficient and precision-adaptive framework.

Similar to the case of the algorithm for computing \(I_\nu(z)\)~\cite{Zaghloul_Johnson_2025}, our present algorithm for computing \(K_\nu(z)\) overcomes the limitations of Algorithm 644, offering improved reliability, an extended domain of computation, and the ability to support quadruple precision where required. Comparison with high-precision values computed using \texttt{Maple} is used to confirm the accuracy of our implementation, and benchmarks show substantial speedup compared to Algorithm 644 in the double precision regime.

This work, similar to part I, contributes to our broader project of developing an open-access, portable, and extensible multiple-accuracy library for the computation of Bessel and other special functions~\cite{Zaghloul2023a,Zaghloul2023b,Zaghloul2024,Zaghloul_Johnson_2025}, tailored to the needs of modern scientific computing.

\section{Algorithm}\label{sec:algorithm}
The modified Bessel function of the second kind, \( K_\nu(z) \), for a complex argument \( z \) and a real order \( \nu \), is computed using different methods depending on the magnitudes of \( \nu \) and \( z \), as well as their relative locations in the complex plane. This section describes the computational approaches and the criteria under which each method is applied to ensure both numerical stability and computational efficiency.

The modified Bessel functions of the second kind satisfy the recurrence relation:
\begin{equation}
K_{\nu+1}(z) = \frac{2\nu}{z} K_\nu(z) + K_{\nu-1}(z)\label{eq:recurrence}.
\end{equation}
This recurrence relation is stable in the forward direction~\cite{Gautschi1967}, but it is numerically unstable in the backward direction. Therefore, it can be reliably used to compute \( K_{\nu+n}(z) \) for \( n = 2, 3, \dots \) when both \( K_\nu(z) \) and \( K_{\nu+1}(z) \) are known.

\begin{figure}[H]
    \centering
    \includegraphics[width=0.6\linewidth]{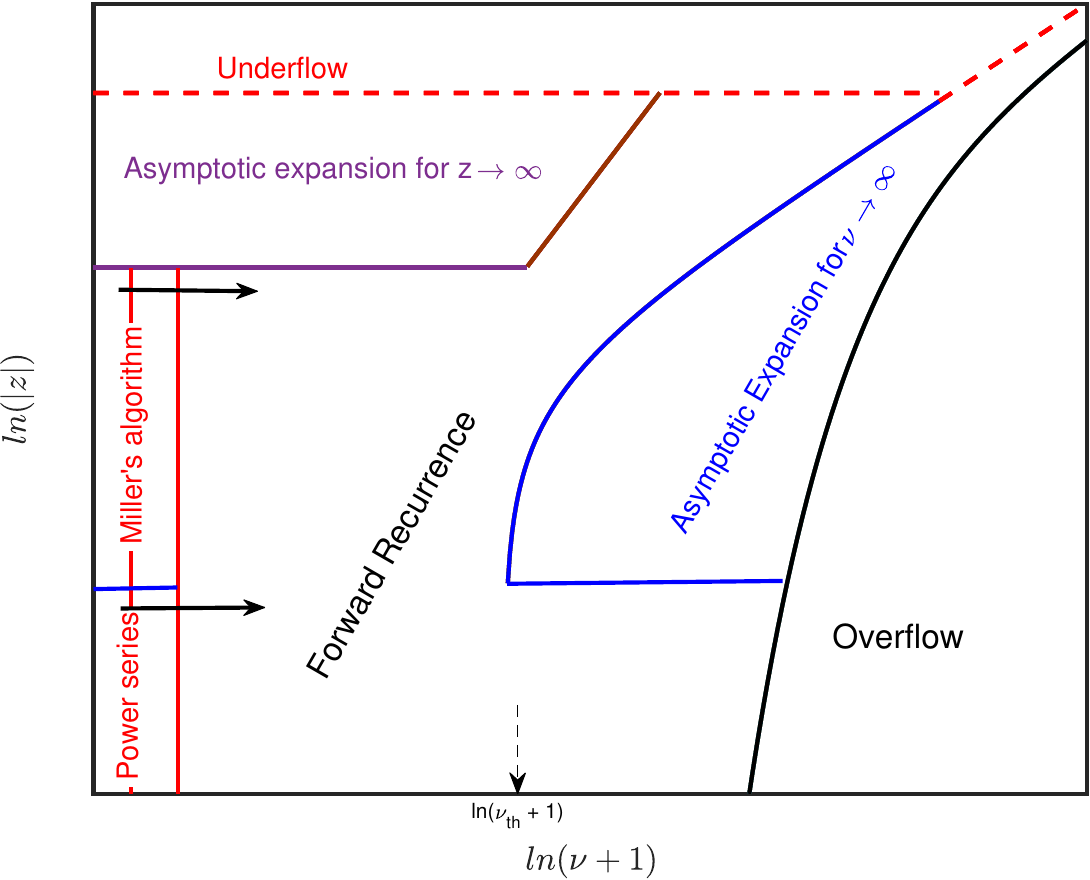} 
    \caption{Schematic computational regions and corresponding methods used in the
present algorithm. }
    \label{fig:Alg_rgns}
\end{figure}

\noindent
Closed-form expressions exist for specific half-integer orders. Such formulas are used herein when \( \nu = \pm \frac{1}{2} \) or \( \nu = \pm \frac{3}{2} \) where,
\begin{equation}
K_{\pm 1/2}(z) = e^{-z} \sqrt{\frac{\pi}{2z}}, \quad
K_{\pm 3/2}(z) = e^{-z} \sqrt{\frac{\pi}{2z}} \left(1 + \frac{1}{z}\right).
\end{equation}

\noindent
Figure 1 provides an overview of the regions and computational methods employed herein in the computation of \( K_\nu(z) \). Pseudo-algorithms for each of the methods used for the regions described below, are provided in the Appendix to facilitate understanding and implementation.

\subsection{Small \( z\): Power Series Expansion \& Forward Recurrence }\label{sec:powerseries}

For small arguments \( z \) and orders \( |\nu| \leq \tfrac{1}{2} \), the convergent series expansion used by Temme~\cite{Temme1975} and Amos~\cite{Amos1983a}, together with the forward recurrence relation, ~\eqref{eq:recurrence}, is sufficient to compute the modified Bessel function of the second kind \( K_\nu(z) \) for \( 0<|z| \leq |z|_{1} \). Amos, in Algorithm 644, uses \( |z|_{1} = 2.0 \) for computations in double precision. However, in this work, we use \( |z|_{1} = 2.2 \) for double precision and \( 3.2 \) for quadruple precision.

\noindent
The convergent series expansion is given by:
\begin{align}
K_\nu(z) &= \sum_{k=0}^\infty C_k f_k, \\
K_{\nu+1}(z) &= \frac{2}{z} \sum_{k=0}^\infty C_k (p_k - k f_k), \label{eq:series},
\end{align}
where
\begin{align}
C_0 &= 1, \qquad C_k = \left( \frac{z^2}{4} \right)^k \frac{1}{k!}, \label{eq:Ck} \\
f_0 &= \frac{\pi \nu}{\sin(\pi \nu)} \left[ \Gamma_1(\nu) \cosh(\mu) + \Gamma_2(\nu) \ln\left( \frac{2}{z} \right) \frac{\sinh(\mu)}{\mu} \right], \label{eq:f0} \\
\mu &= \nu \ln\left( \frac{2}{z} \right), \label{eq:mu}\\
f_k &= \frac{k f_{k-1} + p_{k-1} + q_{k-1}}{k^2 - \nu^2}, \quad k \geq 1. \label{eq:fk}
\end{align}
with the following relations:
\begin{align}
p_0 &= \frac{1}{2} \left( \frac{z}{2} \right)^{-\nu} \Gamma(1+\nu), \quad
p_k = \frac{p_{k-1}}{k - \nu}, \quad k \geq 1, \label{eq:pk} \\
q_0 &= \frac{1}{2} \left( \frac{z}{2} \right)^{\nu} \Gamma(1-\nu), \quad
q_k = \frac{q_{k-1}}{k + \nu}, \quad k \geq 1, \label{eq:qk} \\
\end{align}

The coefficients \( \Gamma_1(\nu) \) and \( \Gamma_2(\nu) \) are defined by:
\begin{align}
\Gamma_1(\nu) &= \frac{1}{2\nu} \left( \frac{1}{\Gamma(1 - \nu)} - \frac{1}{\Gamma(1 + \nu)} \right), \\
\Gamma_2(\nu) &= \frac{1}{2} \left( \frac{1}{\Gamma(1 - \nu)} + \frac{1}{\Gamma(1 + \nu)} \right). \label{eq:gamma12}
\end{align}

As \( \nu \to 0 \), \( \Gamma_1(\nu) \) has the limiting expansion:
\begin{equation}
\Gamma_1(\nu) \sim -\gamma - \left( \frac{\gamma^3}{3} + \gamma \zeta(2) + \frac{2}{3} \zeta(3) \right) \nu^2 + \mathcal{O}(\nu^4), \label{eq:gamma1-limit}
\end{equation}
where \( \gamma \) is the Euler–Mascheroni constant and \( \zeta \) is the Riemann zeta function. For \( \nu^2 \leq \varepsilon \), we use ~\eqref{eq:gamma1-limit} with further \( \nu \to 0 \) simplifications in the calculations. For \( \varepsilon < \nu^2 \leq 10^{-6} \), \( \Gamma_1(\nu) \) is approximated by a Taylor series expansion:
\begin{equation}
\Gamma_1(\nu) = -\sum_{k=0}^N d_{2k+1} \nu^{2k}, \label{eq:gamma1-taylor}
\end{equation}
with precomputed coefficients \( d_{2k+1} \). Twenty-nine such coefficients were generated to 40 significant digits in order to approximate \( \Gamma_1(\nu) \) accurately in the range \( \nu \leq 0.1 \), though only the subrange \( \varepsilon < \nu^2 \leq 10^{-6} \) is used in this implementation.
%
%
%
%
%
\subsection{Uniform Asymptotic Expansion for Large Orders (\(\nu \to \infty\)) \label{sec:largeorder}}
For large orders or for arguments near turning points, the use of asymptotic expansions becomes essential for the accurate evaluation of special functions. For the modified Bessel function of the second kind, \(K_\nu(z)\), the uniform asymptotic expansion provides a reliable and efficient method, especially in the regime of large \(\nu\).

Analogous to the treatment of \(I_\nu(z)\)~\cite{Zaghloul_Johnson_2025}, the asymptotic form of \(K_\nu(\nu z)\) for large \(\nu\) is given by DLMF 10.41.4~\cite{NIST2025}:
\begin{equation}
K_\nu(\nu z) \sim \left( \frac{\pi p}{2 \nu} \right)^{1/2} e^{-\nu \eta}
\sum_{k=0}^{\infty} \frac{(-1)^k U_k(p)}{\nu^k},
\label{eq:nu_inf_expression}
\end{equation}
where the scaling functions \(p\) and \(\eta\) are defined by
\begin{equation}
p = \frac{1}{\sqrt{1+z^2}}, \quad
\eta = \frac{1}{p} + \ln \left( \frac{p z}{1 + p} \right).
\label{eq:p_eta}
\end{equation}

The functions \(U_k(p)\) are polynomials computed recursively to capture higher-order corrections:
\begin{equation}
U_{k+1}(p) = \frac{1}{2} p^2 (1 - p^2) \frac{dU_k}{dp}
+ \frac{1}{8} \int_{0}^{p} \left( 1 - 5 t^2 \right) U_k(t) \, dt, 
\quad U_0(p) = 1.
\label{eq:U_coeff}
\end{equation}
Fortunately, the coefficients required here to generate the polynomials \(U_k(p)\)  coincide with those used in Part~I for computing \(I_{\nu}(z)\); they are precomputed and cached, enabling efficient reuse and improving runtime performance.

This uniform asymptotic expansion offers several advantages: it avoids the instability common in backward recurrence while being more efficient than forward recurrence, and it maintains high accuracy even as \(z\) approaches critical transition points.

The threshold for applying the expansion, \(\nu_{\text{th}}(|z|)\), was determined empirically by comparing the expansion’s results with high-precision reference values generated using \texttt{Maple}~\cite{MAPLE2015}. Although not strictly optimal, this threshold serves as a practical and robust switching criterion, expressed as
\begin{equation}
\nu_{\text{th}}(|z|) = \nu_1 + |z|,
\label{eq: nu_threshold}
\end{equation}
where \(\nu_1\) is a precision-dependent constant. In this study, we use \(\nu_1 = 25.0\) for double precision and \(\nu_1 = 90.0\) for quad precision.

As detailed in~\cite{Zaghloul_Johnson_2025}, the stability of the expansion may initially depend on the phase of \(z\), i.e., \(\arg(z)\), but this dependence vanishes as \(\nu\) becomes sufficiently large. The empirically chosen threshold \(\nu_{\text{th}}(|z|)\) ensures that the expansion remains effectively phase-independent in the large-order limit.
As \( \nu \) increases over the positive real numbers for a fixed nonzero value of \( z \), the pre-sum factor in~\eqref{eq:nu_inf_expression}, namely \( \left( \frac{\pi p}{2\nu} \right)^{1/2} e^{-\nu \eta} \), may overflow. To prevent this, computations are halted and an error message is issued when \( \left| \ln \left( \left( \frac{\pi p}{2\nu} \right)^{1/2} e^{-\nu \eta} \right) \right| > \ln R_{\text{max}} \), where \( R_{\text{max}} \) is the largest representable real number in the arithmetic precision under consideration.
Conversely, as \( z \) increases, the same expression may underflow if
\begin{equation}
\left| \ln \left( \left( \frac{\pi p}{2\nu} \right)^{1/2} e^{-\nu \eta} \right) \right| < \ln R_{\text{min}},
\end{equation} 
where \( R_{\text{min}} \) is the smallest positive normalized real number representable in the given precision. Similarly, for such a case , computations are halted and an error message is issued.

%
%
\subsection{Asymptotic Expansion for Large Argument \( z \to \infty \) and Forward Recurrence}
The asymptotic expansion of the modified Bessel function of the second kind, \( K_{\nu}(z) \), for large argument \( z \to \infty \), is a crucial tool in numerical computation, particularly when \( |z| \gg \nu \). It enables efficient and accurate evaluation of \( K_{\nu}(z) \) in the asymptotic regime defined by
\begin{equation}
|z| \geq \max(|z|_2, \nu^2 / 2).
\label{eq:z_in_border1}
\end{equation}
Here we take \( |z|_2 = 625.0 \) as a precision-independent threshold for applying the large-\(z\) asymptotic expansion (\( z \to \infty \)). This threshold was determined empirically and has some flexibility; it should be regarded as a heuristic rather than a strict boundary.

The asymptotic expansion for large \( z \to \infty \) with fixed \( \nu \) can be expressed as [DLMF 10.40.2] ~\cite{NIST2025},
\begin{equation}
K_\nu(z) \sim \left( \frac{\pi}{2z} \right)^{1/2} e^{-z} \sum_{k=0}^\infty \frac{a_k(\nu)}{z^k}, \qquad |\arg z| \leq \tfrac{3}{2}\pi - \delta
\label{eq:asymptotic1}
\end{equation}
where
\begin{equation}
a_k(\nu) = \frac{(4\nu^2 - 1^2)(4\nu^2 - 3^2) \cdots (4\nu^2 - (2k - 1)^2)}{k! \, 8^k}
\label{eq:ak}
\end{equation}

~\eqref{eq:asymptotic1} can be calculated recursively and written as:
\begin{equation}
K_\nu(z) \sim \left( \frac{\pi}{2z} \right)^{1/2} e^{-z} \sum_{k=0}^\infty T_k, \qquad |\arg z| \leq \tfrac{3}{2}\pi - \delta
\label{eq:asymptotic2}
\end{equation}
with the recurrence:
\begin{align}
T_0 &= 1, \\
T_{k+1} &= \frac{4\nu^2 - (2k+1)^2}{8(k+1)z} \, T_k, \qquad k \geq 0
\label{eq:Tk_recursion}
\end{align}

A sufficient condition for the series to converge is:
\begin{equation}
\left| \frac{T_{k+1}}{T_k} \right| = \left| \frac{4\nu^2 - (2k+1)^2}{8(k+1)z} \right| \leq 1
\label{eq:ratio_condition}
\end{equation}

This ratio decreases as \( k \) increases, and the condition for convergence for large \( z \) is \(|z| > \nu^2 / 2.0\). 

The summation in \eqref{eq:asymptotic2} is terminated after a prescribed maximum number of terms, following the same strategy used in Part~I for \(I_{\nu}(z)\).

For the region of applicability of the asymptotic expression given in ~\eqref{eq:asymptotic1}, underflow may occur when:
\begin{equation}
\Re \left[ \ln\left( \left( \frac{\pi}{2z} \right)^{1/2} e^{-z} \right) \right] < \ln(R_{\min})
\label{eq:underflow}
\end{equation}

\noindent
As evident from the condition in ~\eqref{eq:underflow}, underflow in this region of the computational domain depends on the angle (phase) of the complex variable \( z \). Thus, there is no well-defined, single-valued boundary in the \( \nu \)-\(|z| \) plane. Consequently, the dashed line used for the upper limit of this region in Fig.~1 is symbolic rather than a strict mathematical boundary.

For the special case of a real argument, \( z = x \), the inequality in ~\eqref{eq:underflow} is solved iteratively, showing that underflow occurs when
$x > 705.34 \text{ (for double precision)},\ x > 11352.081 \text{ (for quad precision)}$

\noindent
Accordingly, the region of applicability of this asymptotic expansion for large argument \( z \to \infty \) with fixed \( \nu \), for both double and quad precision,  is therefore given by:
\begin{equation}
|z| > \max\left( \frac{\nu^2}{2.0}, \, 625.0 \right)
\label{eq:region_validity}
\end{equation}
provided that the underflow condition in ~\eqref{eq:underflow} is not reached.
It is worth noting that, according to Amos~\cite{Amos1983a}, the method used for computing \( K_\nu(z) \) in the region \( |z| > 2.0 \) is Miller’s algorithm combined with forward recurrence (explained in the following subsection). However, in our implementation, we employed the asymptotic expansion for \( z \to \infty \) in the region defined by ~\eqref{eq:region_validity}, as it is not only more efficient but also readily extendable to a wider computational domain using quad-precision arithmetic.
For the intermediate region where \( 625.0 < |z| \leq \nu^2 / 2.0 \) and \( \nu < \nu_{\text{th}}(|z|) \), forward recurrence using ~\eqref{eq:recurrence}, which is numerically stable, is employed to compute \( K_\nu(z) \).
%
%
\subsection{Intermediate Region} \label{sec:intermediate}

In the intermediate region of \( |z| \), specifically for \( |z|_1 < |z| < |z|_2 \) with \( \nu \leq \nu_{\text{th}}(|z|) \), and using the notation of Abramowitz and Stegun~\cite{Abramowitz1964}, the modified Bessel function of the second kind \( K_\nu \) can be written as
\begin{equation}
K_\nu(z) = \sqrt{\pi} \, (2z)^\nu \, e^{-z} \, U\left(\nu+\tfrac{1}{2}, 2\nu+1, 2z\right),
\label{eq:Knu_as_U}
\end{equation}
where \( U(a,b,z) \) is a confluent hypergeometric function, which for \( \Re(z) > 0 \) and \( \Re(a) > 0 \) may be defined as
\begin{equation}
\Gamma(a) \, U(a,b,z) = \int_0^\infty e^{-zt} \, t^{a-1} (1+t)^{b-a-1} \, dt.
\label{eq:U_integral}
\end{equation}

The function \(k_n(z)\) is related to \(U(a,b,z)\) and hence to \(K_\nu(z)\) and can be expressed as 
\begin{equation}
k_n(z) = (-1)^n \left[ \frac{1}{\pi n!} (-1)^n \cos(\nu\pi) \Gamma\left(\tfrac{1}{2} + \nu + n\right) \Gamma\left(\tfrac{1}{2} - \nu + n\right) \right] U\left(\nu + \tfrac{1}{2}, 2\nu + 1, 2z\right)
\label{eq:kn_def}
\end{equation}
The function \(k_n(z)\) satisfies the recurrence relation ~\cite{Temme1975,Campbell_81}
\begin{equation}
k_{n+1}(z) - b_n \, k_n(z) + a_n \, k_{n-1}(z) = 0, \quad n \geq 1,
\label{eq:kn_recurrence}
\end{equation}
with
\begin{equation}
a_n = \frac{(n - \tfrac{1}{2})^2 - \nu^2}{n^2 + n}, \quad
b_n = \frac{2(n + z)}{n+1}, \quad n = 1,2,\ldots
\label{eq:an_bn}
\end{equation}

For large values of \( n \), the function \( k_n(z) \) behaves as
\begin{equation}
k_n(z) \sim \pi^{-1/2} \cos(\pi\nu) \, 2^{1/4} \, n^{-1/2} \, z^{-\nu - 1/4} \exp\left[z - 2(2nz)^{1/2}\right].
\label{eq:kn_asymptotic}
\end{equation}

Hence, Miller’s algorithm with backward recurrence can be used to determine \( k_0(z) \) and \( k_1(z) \) with the help of a normalization relation. In Miller’s algorithm, a terminal value \( n = M \) is chosen such that one sets
\begin{equation}
k_{M+1}^M = 0, \quad k_M^M = \varepsilon,
\label{eq:limits_Miller}
\end{equation}
and applies backward recurrence along with the normalization condition:
\begin{equation}
\sum_{n=0}^{M} k_n^M(z) = (2z)^{-\nu - 1/2} = S_M,
\label{eq:normalization}
\end{equation}
to find \( k_0(z) \) and \( k_1(z) \). Having evaluated \( k_0(z) \) and \( k_1(z) \), the values of \( K_\nu(z) \) and \( K_{\nu+1}(z) \) can be computed using
\begin{equation}
K_\nu(z) \approx \sqrt{\frac{\pi}{2z}} e^{-z} \frac{k_0^M(z)}{S_M},
\label{eq:Knu_from_k0}
\end{equation}
and
\begin{equation}
K_{\nu+1}(z) = K_\nu(z) \left[\nu + z + \tfrac{1}{2} - \frac{k_1^M(z)}{k_0^M(z)}\right] / z.
\label{eq:Knu_plus1}
\end{equation}

Based on extensive numerical studies, Amos~\cite{Amos1983a} derived the following expression for determining a suitable value of \( M \) for small \( |z| \),
\begin{equation}
M = \frac{0.97}{2|z|} \times \left( \frac{\ln C + \dfrac{|z| \cos(A\theta)}{1 + 0.008|z|}}{2 \cos(B\theta)} \right)^2 + 1.5,
\label{eq:M_estimate}
\end{equation}
where
\begin{equation}
\theta = \arg z, \quad
A = \frac{3}{1 + |z|}, \quad
B = \frac{0.49 \cdot 30}{28 + |z|}, \quad
C = \frac{4 \cos(\pi\nu)}{\varepsilon \pi^{1/2} (2|z|)^{1/4}}.
\label{eq:ABC}
\end{equation}

\noindent
In contrast to Amos' implementation, the following modifications are proposed in the present work to improve computational efficiency.

\noindent
First, the term $\ln C$ varies only slowly over the region of interest and has a relatively weak influence on the value of $M$, which itself serves as a heuristic truncation parameter in Miller's algorithm. Consequently, its direct evaluation can be avoided by replacing it with a fixed constant determined empirically. Extensive numerical testing indicates that using $\ln C = 32$ for double precision and $\ln C = 56$ for quadruple precision provides a stable and reliable approximation that preserves the desired accuracy while improving computational efficiency. In this sense, the constant approximation can be viewed as a practical lower-bound estimate consistent with the prescribed accuracy requirement.

This simplification is applied specifically to $\ln(C)$ due to its weak sensitivity. Other terms in ~\eqref{eq:M_estimate} exhibit stronger dependence on $|z|$ and $\nu$, and therefore must be retained in their original form to ensure the accuracy and robustness of the method.

\noindent
Second, the expression in~\eqref{eq:M_estimate} is used herein only within the region $|z| \le |z|_{3}$, for which our tests indicate that
$|z|_{3} = 40$ in double precision and $|z|_{3} = 160$ in quadruple precision are appropriate choices.
For larger arguments, $|z| > |z|_{3}$, we instead employ fixed values of $M$, namely
$M = 7$ for double precision and $M = 16$ for quadruple precision.

The approaches discussed in Sections \ref{sec:powerseries} and \ref{sec:intermediate} share some similarities with \cite{Temme1975,Campbell_81,Amos1986}, particularly in terms of their general computational framework. However, the present algorithm introduces a few effective modifications similar to the ones described above in addition to the selection of domain boundaries to improve computational efficiency. Additionally, our \texttt{Fortran} implementation incorporates several optimization techniques; for example, the hyperbolic functions of complex arguments are avoided by computing the exponential function once and reusing it where appropriate.

Furthermore, the direct use of asymptotic expansions for large values of $\nu$ and $z$, while not explicitly implemented in earlier algorithms such as those by Temme~\cite{Temme1975}, Campbell~\cite{Campbell_81}, Thomson and Barnett~\cite{Thompson_Barnett_1987}, Amos~\cite{Amos1983a, Amos1986}, or Press and Teukolsky~\cite{Press_Teukolsky_1991}, has proven particularly beneficial in our work. This strategy facilitated accurate and efficient extension of the algorithm over a broader computational domain, especially in the context of quadruple-precision arithmetic.

\subsection{Analytic continuation to the left half–plane \texorpdfstring{$(\Re z<0)$}{(Re z<0)}}
\label{subsec:cont-left-half}

On the principal branch, \(K_\nu(z)\) has a branch cut along the negative real axis.  
Across this cut, the standard connection formula (DLMF 10.34.2) ~\cite{NIST2025}, gives for \(z\) off the cut,
\begin{equation}
\label{eq:DLMF-10342}
K_\nu\!\big(z\,e^{\pm i\pi}\big) = e^{\mp i\pi \nu}\,K_\nu(z)\;\mp\; i\pi\,I_\nu(z).
\end{equation}

For \(\Re z<0\), we set \(w=-z\) so that \(\Re w>0\). Approaching \(z\) from the upper half–plane (\(\Im z>0\)) corresponds to rotating \(w\) by \(+\pi\), while approaching from the lower half–plane (\(\Im z<0\)) corresponds to rotating \(w\) by \(-\pi\). Thus,
\begin{equation}
z = w\,e^{i s \pi}, \qquad s=\operatorname{sgn}(\Im z)\in\{+1,-1\}.
\label{eq:z_w_mapping}
\end{equation}

Applying \eqref{eq:DLMF-10342} with \(z\mapsto w\) and the sign determined by \(s\) as in \eqref{eq:z_w_mapping} yields the practical continuation (valid off the branch cut):
\begin{equation}
K_\nu(z)
= e^{-i s \pi \nu}\,K_\nu(-z)\;-\; i\,s\,\pi\,I_\nu(-z),
\qquad s=\operatorname{sgn}(\Im z),\ \Re z<0.
\label{eq:K-left-half-plane}
\end{equation}

On the branch cut, we write \(z=-x\) with \(x>0\). Then
\begin{equation}
\label{eq:cut-limits}
K_\nu(-x \pm i0)
=
e^{\mp i\pi \nu}\,K_\nu(x)\;\mp\; i\pi\,I_\nu(x), 
\qquad x>0,
\end{equation}
so an implementation must consistently choose the upper bank (\(+\)) or lower bank (\(-\)) when \(\Im z=0\). 

In addition, direct evaluation of \(e^{-i s\pi \nu}\) in \eqref{eq:K-left-half-plane} can lose accuracy for large \(|\nu|\) or near integer \(\nu\). Accordingly for numerically stable phase, we  write \(\nu=n+f\) with \(n=\lfloor \nu\rfloor\in\mathbb{Z}\) and \(f\in[0,1)\); then
\begin{equation}
\label{eq:phase-stable}
e^{-i s\pi \nu}
=
(-1)^n\!\big[\cos(\pi f)\;-\; i\,s\,\sin(\pi f)\big],
\end{equation}
which is well conditioned, including as \(f\to 0\) (integer orders).  
In particular, for integer \(n\),
\begin{equation}
\label{eq:int-order}
K_n(z)=(-1)^n K_n(-z)\;-\; i\,s\,\pi\,I_n(-z),\qquad \Re z<0,\ s=\operatorname{sgn}(\Im z).
\end{equation}


\section{Implementation and Accuracy Verification}\label{sec:implementation}

The modular design of the \texttt{Fortran} implementation of the algorithm from Section~\ref{sec:algorithm} enables support for both double and quadruple precision arithmetic. Control over the precision is provided by the integer parameter \texttt{rk}, which defines the real kind used throughout the module and its driver. This parameter is specified within a dedicated module named \texttt{set\_rk}.

A comprehensive evaluation was carried out to assess both the accuracy and performance of the proposed algorithm for computing modified Bessel functions of the second kind. This involved detailed comparisons with Algorithm~644, a well-established\texttt{ Fortran} implementation known for its double-precision capabilities. To the best of our knowledge, free or open-source implementations of complex-argument Bessel functions in compiled languages (including \texttt{Fortran}) that provide native support for quadruple-precision arithmetic remain limited, with many existing approaches relying on comparatively slow arbitrary-precision libraries.

To ensure a reliable foundation for accuracy evaluation, high-precision reference values were generated using \texttt{Maple}'s arbitrary-precision capabilities and independently verified against results from the \texttt{mpmath} library~\cite{mpmath}. The relative errors between these reference values and those computed by the present proposed algorithm and Algorithm~644 were systematically analyzed across a wide range of parameters, including both small and large orders and arguments. This thorough assessment demonstrates the robustness and accuracy of the new algorithm under realistic computational scenarios. Furthermore, performance benchmarks reveal that the proposed implementation significantly outperforms Algorithm~644 in terms of efficiency.\\

As discussed in detail in Part~I~\cite{Zaghloul_Johnson_2025}, extending numerical 
algorithms from double to quadruple precision significantly increases the 
available dynamic range and reduces round-off error, but does not eliminate 
all numerical challenges. In particular, exponential growth in asymptotic 
representations can still lead to overflow/underflow, albeit at larger values 
of $|z|$ and $\nu$, and cancellation effects may arise for complex arguments.
\noindent
The same considerations apply to the present algorithm for $K_\nu(z)$. In 
addition, the stability of recurrence relations and the accurate evaluation 
of asymptotic expansions remain critical over the extended parameter domain. 
These challenges are addressed through adaptive region selection, careful 
initialization of recurrence procedures, and controlled truncation strategies, 
ensuring stable and accurate computations in both double and quadruple precision.

\subsection{Accuracy Verification: Double Precision} \label{sec:Acc_Ver_ Dbl_Precision}

\begin{figure}[H]
    \centering
    \includegraphics[width=0.6\linewidth]{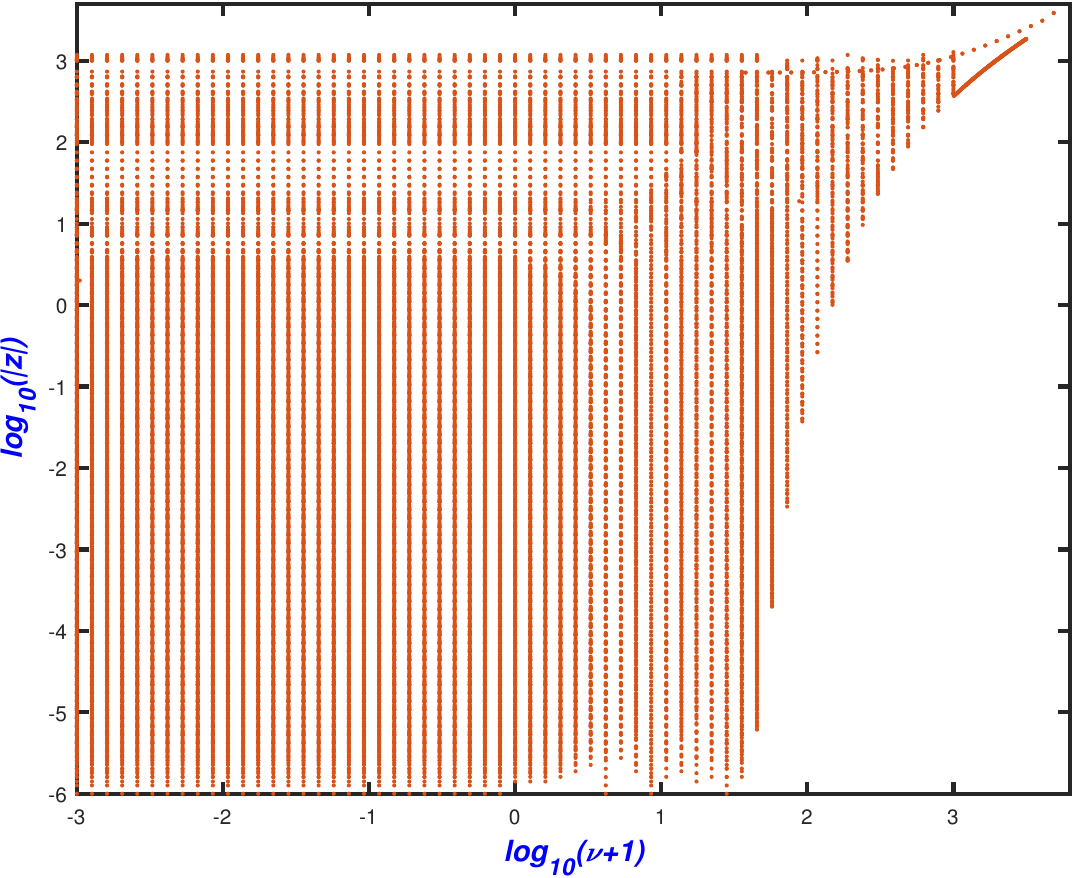} 
    \caption{Distribution of the test grid computed using double-precision arithmetic. For all input values shown, \texttt{Maple} evaluations remain within the representable range of real numbers in IEEE double-precision arithmetic.}
    \label{fig:dbl_pts_new}
\end{figure}

To rigorously evaluate the algorithm’s accuracy in double precision, we constructed a comprehensive test grid comprising 283180 points, predominantly distributed quasi-uniformly on a logarithmic scale in the \( \nu \)-\(|z|\) plane. The complete test grid is visualized in Fig. \ref{fig:dbl_pts_new}.

All function values at these test points were computed using high-precision\texttt{ Maple} calculations, ensuring that they fall within the range of representable floating-point numbers in double-precision arithmetic.

\begin{figure}[H]
    \centering
    \subfloat[]{\includegraphics[width=0.48\linewidth]{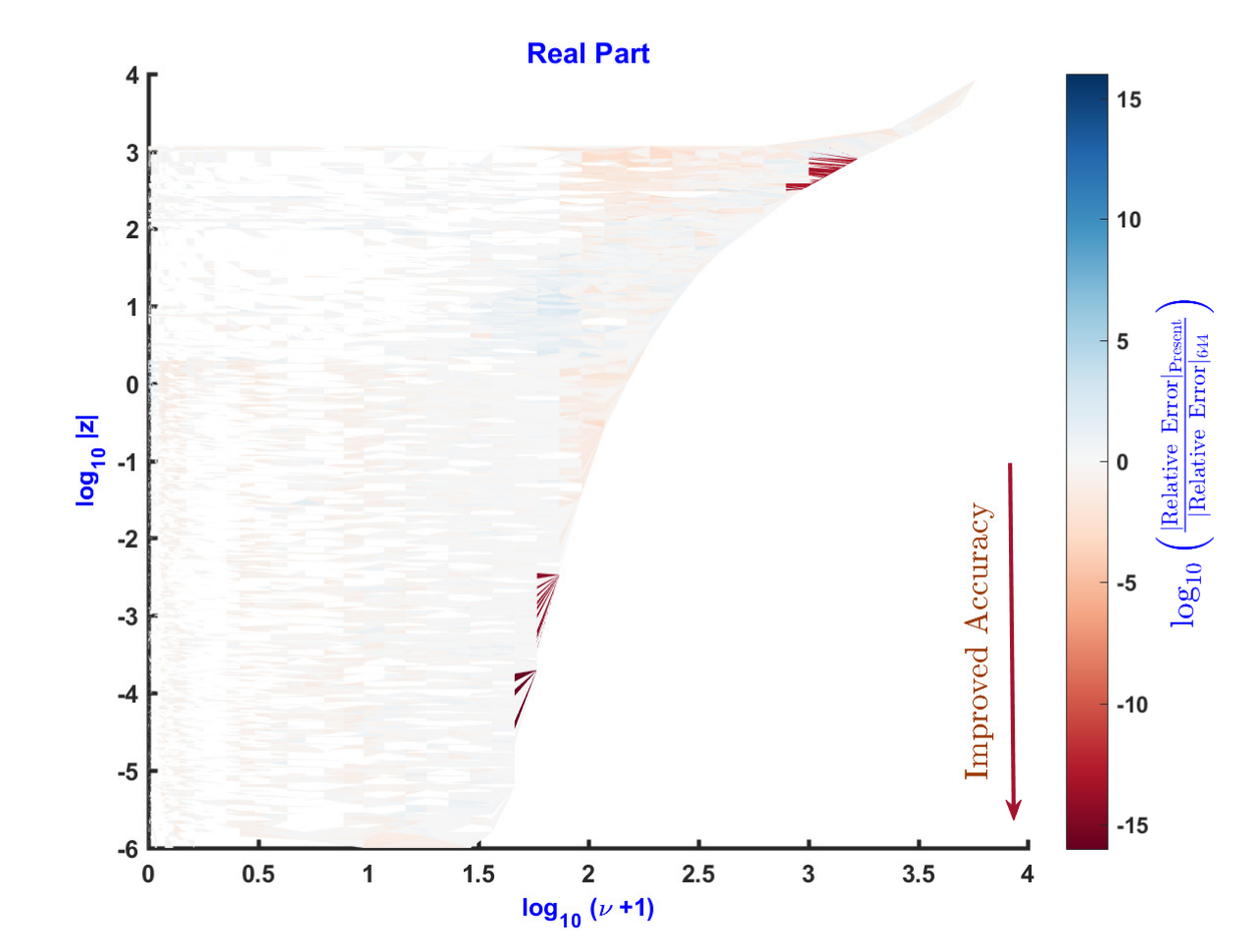}\label{fig:3a}}
    \subfloat[]{\includegraphics[width=0.48\linewidth]{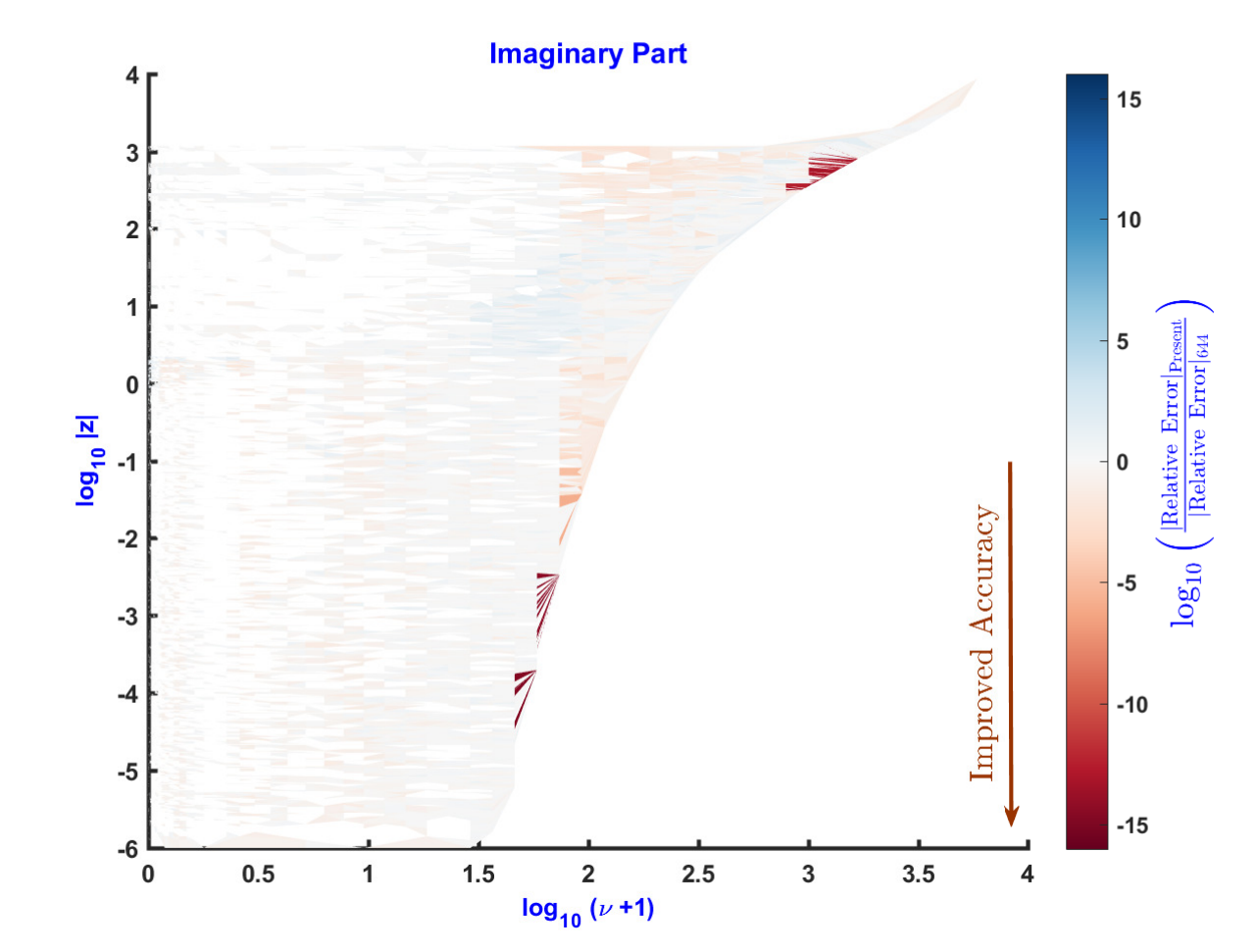}\label{fig:3b}}
    \caption{Colormaps showing numerical accuracy comparison for \(K_{\nu}(z)\) in terms of \(\log_{10}\) of the ratio between the relative error of the present algorithm and that of Algorithm~644, with high-precision \texttt{Maple} values as reference. Results are given for the real part (a) and imaginary part (b)}\label{fig:dp_acc_comp}
\end{figure}

To further illustrate the numerical error statistics, Figure~\ref{fig:dp_acc_comp} shows a visual comparison of the accuracy of the present algorithm and Algorithm~644 over the entire test domain. The figure presents two-dimensional colormap plots of the base-10 logarithm of the \textit{ratio} of the relative errors, \(
\boldsymbol{|\textit{computed} - \textit{reference}| / |\textit{reference}|}
\),
obtained with the present method and with Algorithm~644, using \texttt{Maple}’s arbitrary-precision results as the reference standard. Results are shown separately for the real part~(a) and imaginary part~(b) of \( K_{\nu}(z) \). Positive values in the plots correspond to regions where the present algorithm exhibits larger relative error than Algorithm~644, whereas negative values highlight regions of improved accuracy. The color distributions reveal that, across most of the tested parameter space, the proposed implementation achieves equal or better accuracy, with the largest gains observed in parameter regions that are numerically challenging for double-precision computations.

Figure~~\ref{fig:dp_acc_comp} presents colormap plots of the base-10 logarithm of the ratio of the relative error of the present algorithm to that of Algorithm~644, using the same \texttt{Maple} reference values. Since this quantity may be either positive or 
negative, the color scale is chosen symmetrically about zero. Negative values indicate regions where the present algorithm is more accurate, positive values indicate regions where Algorithm~644 is more accurate, and values near zero indicate comparable performance.
In large portions of the computational domain, Algorithm~644 halts execution and returns an error code because of overly conservative overflow thresholds. These failure regions appear several orders of magnitude before the true overflow boundary. In such cases, Algorithm~644 outputs \((0.0, 0.0)\), leading to a \(100\%\) relative error. This behavior is visible as red areas (corresponding to extremely small relative error ratios) in parts~(a) and~(b) of Fig.~\ref{fig:dp_acc_comp}. By contrast, the present algorithm continues to deliver accurate results even within these challenging regions. 

Away from these regions where Algorithm~644 fails, the accuracy of both algorithms is comparable with the present algorithm producing more accurate results for large \(\nu\).

\subsection{Accuracy Verification: Quad Precision}

The accuracy and domain restrictions discussed in Section~\ref{sec:Acc_Ver_ Dbl_Precision} originate primarily from the finite resolution of double-precision arithmetic (about 16 significant digits). Such limitations become particularly severe for large orders and/or arguments, where numerical instabilities can dominate. While arbitrary-precision arithmetic, as provided by tools like \texttt{Maple} or \texttt{mpmath}, can, in principle, eliminate these issues, its high computational cost makes it impractical for large-scale or time-critical applications. By employing quad-precision arithmetic within a compiled language such as \texttt{Fortran}, the present implementation achieves a substantial extension of the computational domain, maintains numerical stability, and delivers high accuracy with far superior performance to arbitrary-precision approaches. This capability underpins the following assessment of quad-precision accuracy.

\begin{figure}[H]
    \centering
    \includegraphics[width=0.6\linewidth]{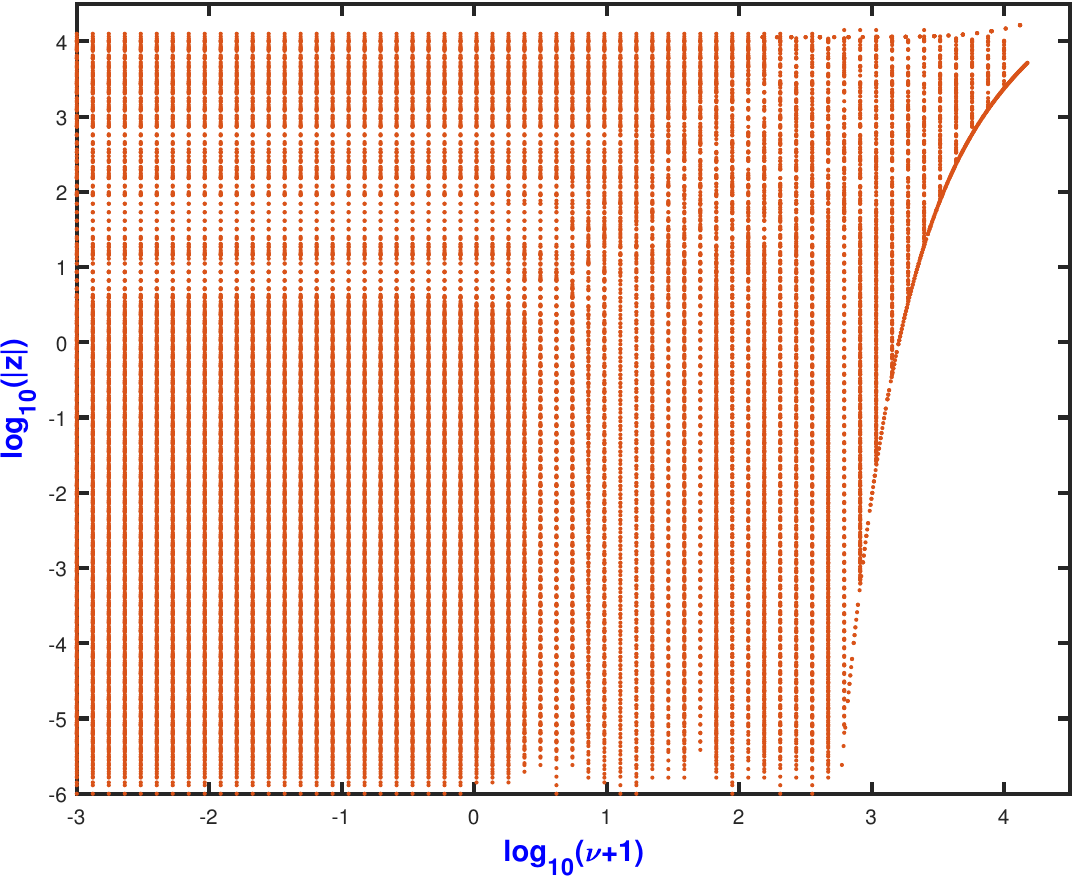} 
    \caption{The grid of tested points using quad-precision arithmetic as shown in Fig.~\ref{fig:Alg_rgns}. \texttt{Maple} calculations for these input points fall within the range of the minimum and maximum floating-point real numbers in quad-precision arithmetic. }
    \label{fig:qp_pts_new}
\end{figure}

To assess the accuracy of the quad-precision implementation, we constructed a new test grid spanning the enlarged \(\nu\)-\(|z|\) domain accessible at this precision level. The grid includes a total of 280{,}496 test points (see Fig.~\ref{fig:qp_pts_new}). Relative to the domain examined in double precision, the quad-precision range extends by more than an order of magnitude in both order and argument.

\begin{figure}[H]
    \centering
    \subfloat[]{\includegraphics[width=0.48\linewidth]{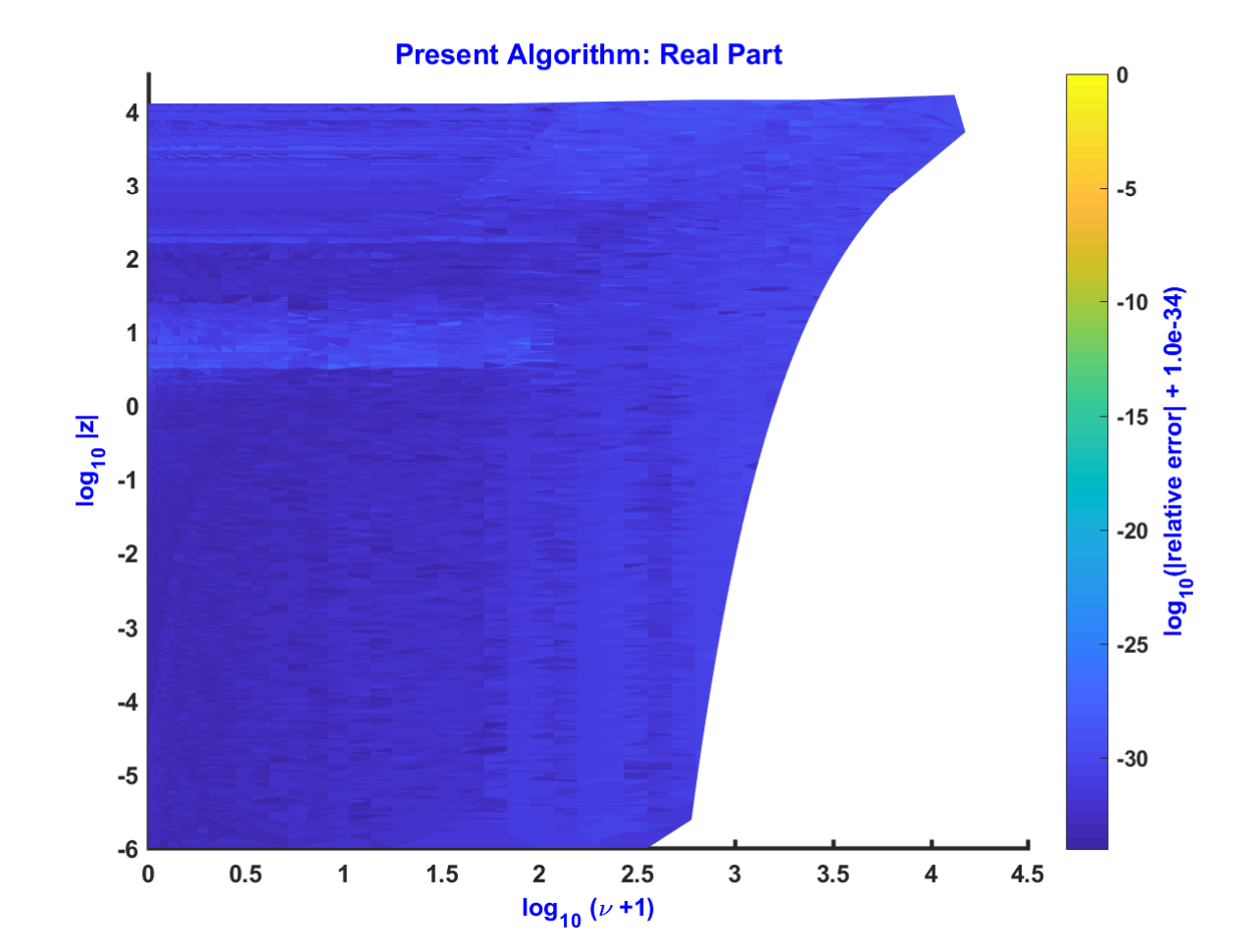}\label{fig:K_qb_acc_real}}
    \subfloat[]{\includegraphics[width=0.48\linewidth]{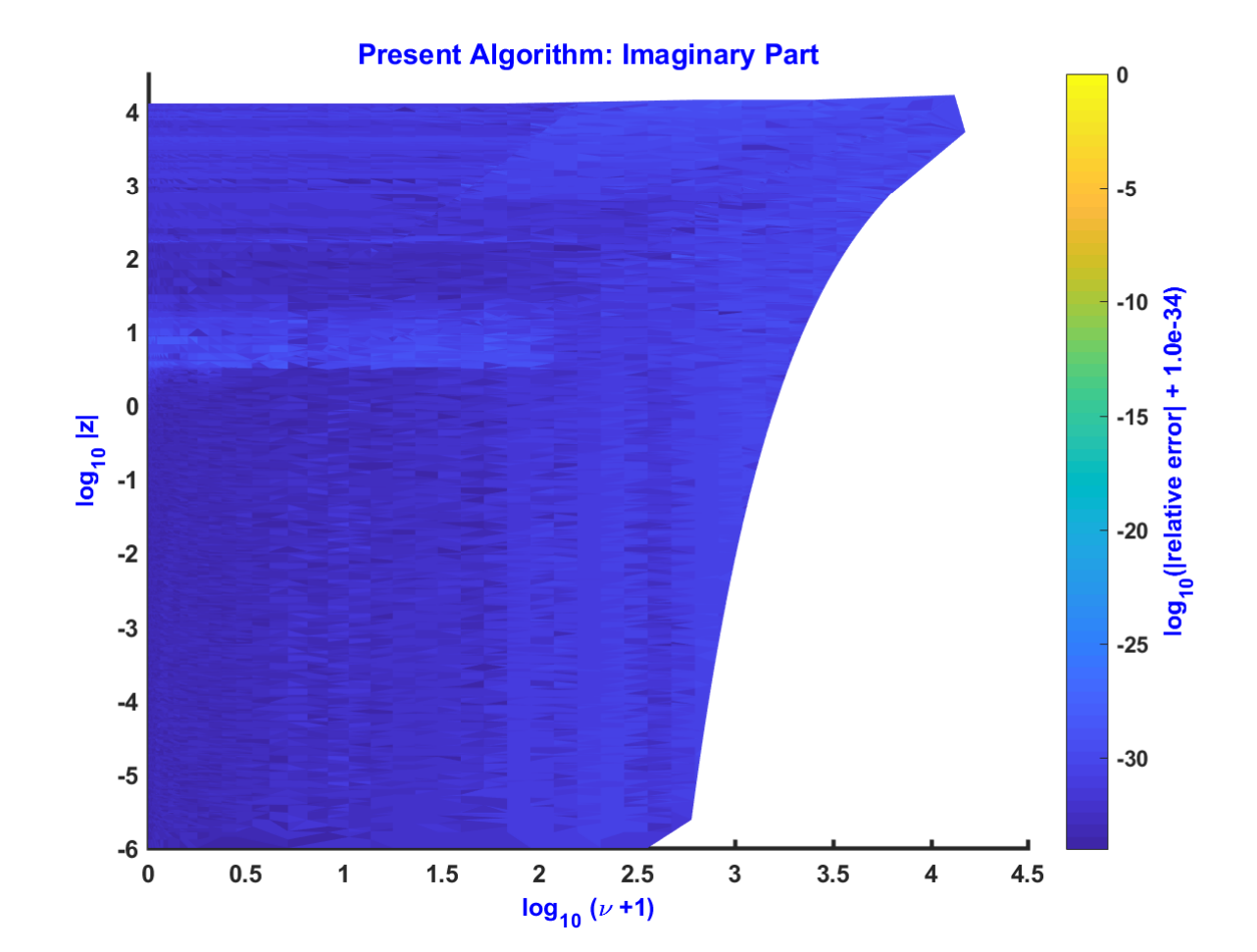}\label{fig:K_qp_acc_Im}}
    \caption{Colormap plots of the relative error in calculating the real part (a) and imaginary part (b) of \( K_{\nu} (z) \) for the dataset points tested using the present algorithm, with \texttt{Maple} calculations as the reference.}
     \label{fig:qp_acc}
\end{figure}

The corresponding accuracy results are summarized in Fig.~\ref{fig:qp_acc}, which displays two-dimensional colormaps of the relative error for the real part~(a) and imaginary part~(b) of \(K_{\nu}(z)\) over the expanded test grid. As before, high-precision \texttt{Maple} calculations serve as the reference standard. The results demonstrate that the present \texttt{Fortran} implementation maintains exceptional accuracy across the domain, with relative errors  \(<10^{-26}\) for the real part and \(<10^{-28}\) for the imaginary part, even in the vicinity of region boundaries.

\section{Efficiency Benchmarking: Present Algorithm vs Algorithm~644}

To evaluate the computational efficiency of the present algorithm for \( K_{\nu}(z) \), a systematic benchmarking study is performed against Algorithm~644. The benchmarking employs the dataset described in Section~3.1, ensuring a comprehensive assessment across the entire computational domain.  

Execution times are measured for both algorithms under identical computational conditions, enabling a fair, point-by-point comparison. Since Algorithm~644 was designed for double precision and does not support quad precision, the efficiency tests of the present algorithm are restricted to double precision. This evaluation highlights the computational advantages of the new algorithm, particularly in terms of execution speed. Support for higher precision, and robustness in regions of extreme parameter values where Algorithm~644 fails have been discussed above.  

Efficiency tests are based on the same set of points used in the double-precision accuracy study. Each test computes all points 50 times and is repeated 21 times. To minimize noise from system-level fluctuations, the shortest execution time out of the 21 repetitions is recorded~\cite{BenchmarkTools}. Timing is measured using the intrinsic \texttt{SYSTEM\_CLOCK} function, which in GNU~\texttt{Fortran} maps to the high-resolution \texttt{QueryPerformanceCounter} system call, providing sufficient resolution to resolve sub-microsecond events.  

Benchmarking experiments are conducted under controlled conditions on an Intel\textsuperscript{\textregistered} Core\textsuperscript{TM} i7-6600U CPU @ 2.60–2.81~GHz with minimal background processes. The following Fortran compilers are tested:  
\begin{itemize}
    \item GNU Fortran (i686-posix-dwarf-rev0, Built by MinGW-W64 project) 8.1.0,
    \item GNU Fortran (Rev3, Built by MSYS2 Project) 12.1.0,
    \item NAG Fortran Compiler Release 7.1 (Hanzomon) Build 7110,
    \item Intel(R) Fortran Intel(R) 64 Version 2021.9.0 Build 20230302\_000000,
    \item \texttt{IFX} (LLVM-based), from Intel(R) Fortran 64 Version 2021.9.0.
\end{itemize}

All tests are compiled with optimization levels \texttt{-O3}, \texttt{-O2}, \texttt{-O1}, and \texttt{-O0} to assess the impact of compiler optimizations. Results are summarized in Table~\ref{tab:execution_time_comparison}.  

The data in Table~\ref{tab:execution_time_comparison} show that the present algorithm consistently outperforms Algorithm~644 across the entire computational domain and within all sub-regions. Depending on the compiler and optimization level, the execution time of the present algorithm ranges between about 49\% and 70\% of that required by Algorithm~644. 

\begin{figure}[H]
    \centering
    \subfloat[]{\includegraphics[width=0.47\linewidth]{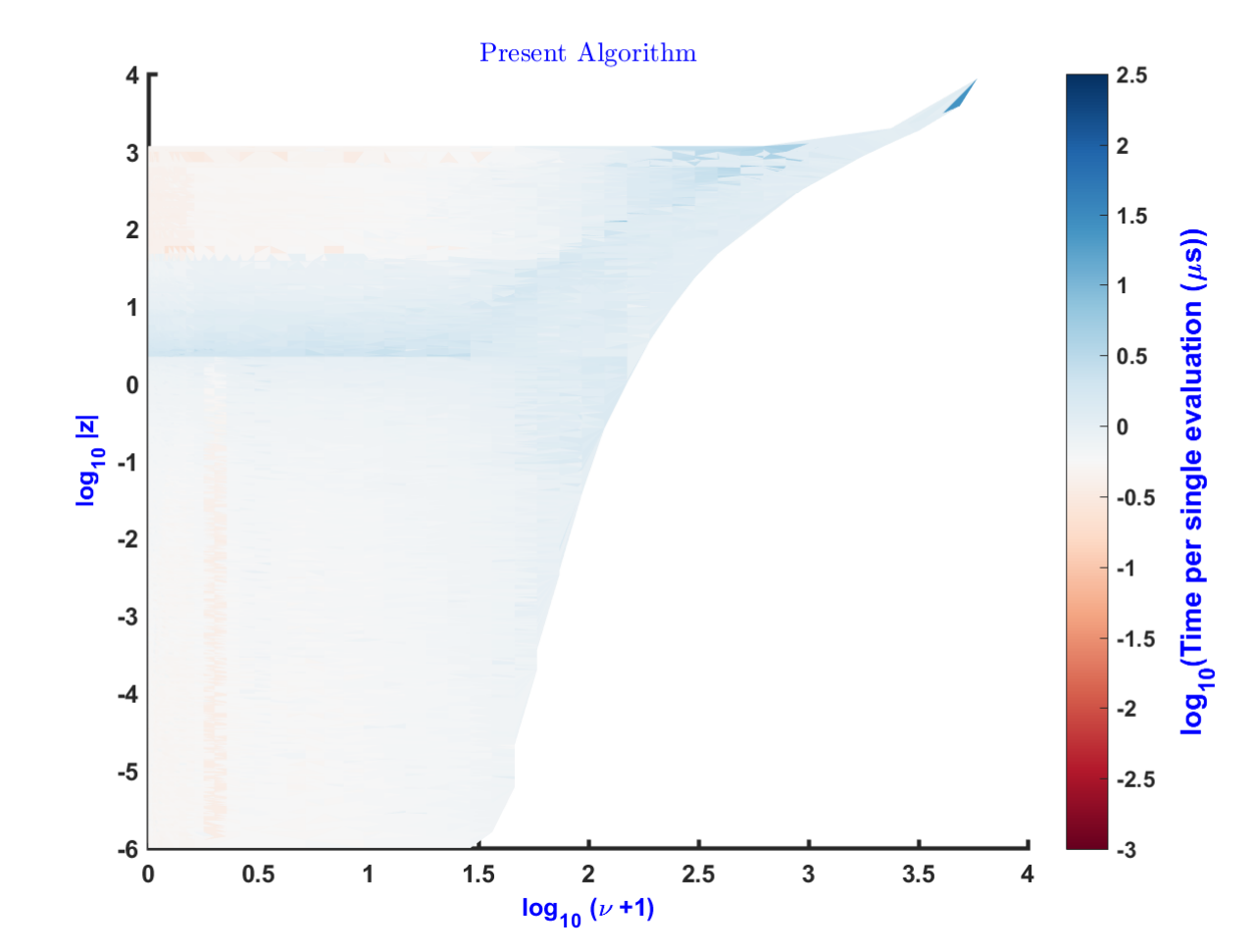}\label{fig:pt_by_pt_a}}
    \subfloat[]{\includegraphics[width=0.47\linewidth]{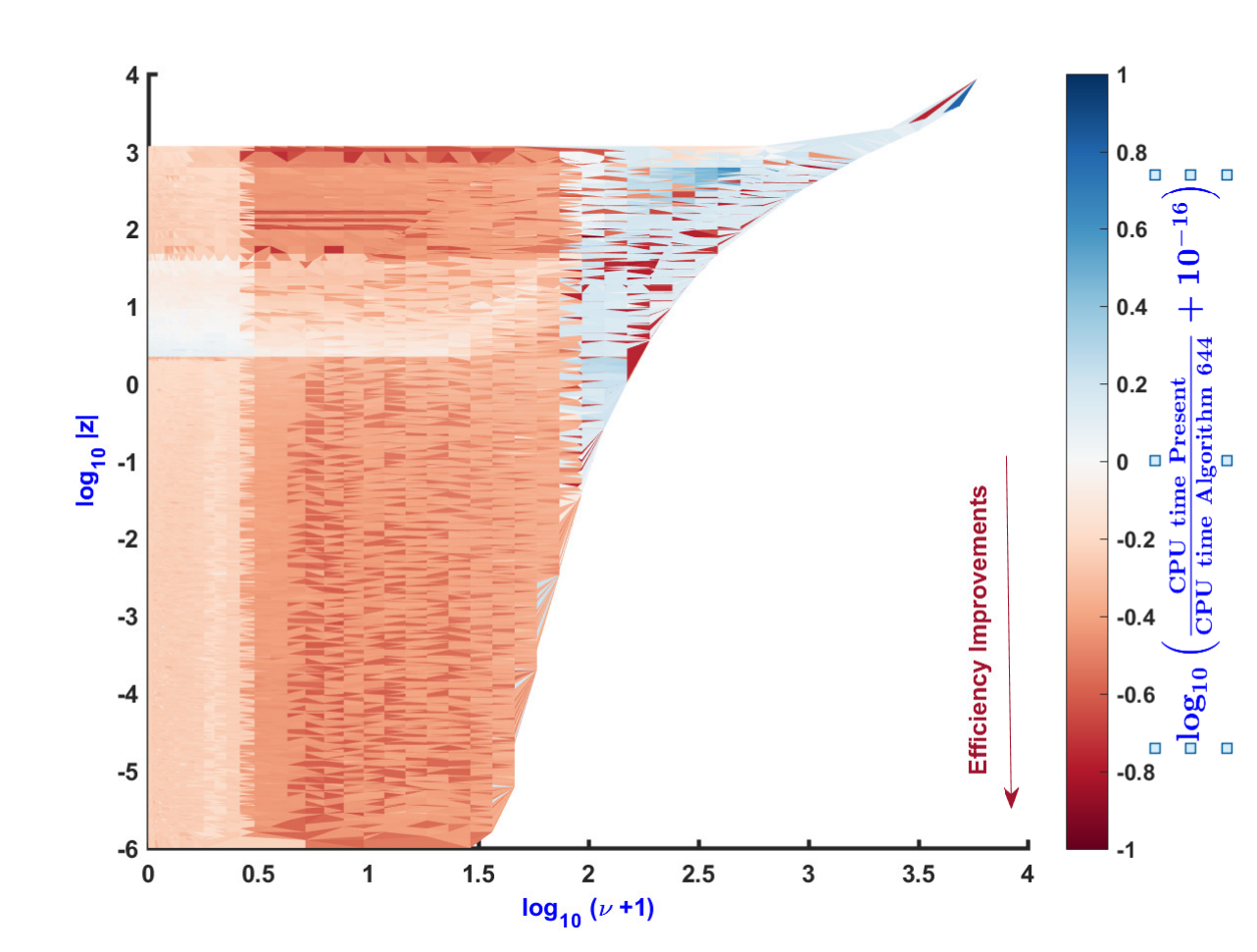}\label{fig:pt_by_pt_b}}
 \caption{Two-dimensional colormap plots of the base-10 logarithm of: \textbf{(a)} execution time per evaluation (in~\(\mu\)s) of the present algorithm, and \textbf{(b)} the time \textit{ratio} of the present algorithm to Algorithm~644 across the domain. Computations used an Intel\textsuperscript{\textregistered} Core\textsuperscript{TM} i7-6600U CPU @ 2.60--2.81~GHz with GNU Fortran 12.1.0 (Rev3, MSYS2) with optimization level \texttt{-O3}.}
\end{figure}

Further insight is obtained by examining execution time at the level of individual test points. Figure \ref{fig:pt_by_pt_a} presents a colormap of \(\log_{10}\) of the execution time (in~\(\mu\)s) for a single evaluation, while Fig. \ref{fig:pt_by_pt_b} displays \(\log_{10}\) of the execution-time \textit{ratio} between the present algorithm and Algorithm~644. The plots show that the present algorithm is faster across most of the domain, with particularly pronounced gains for small–to–intermediate \(\nu\). In the large-order regime (\(\nu \to \infty\)), where the uniform asymptotic expansion is employed together with forward recurrence in the large-\text{z} region, there are isolated regions in which the present algorithm is substantially, often orders of magnitude, faster than Algorithm~644.

\begin{table}[ht]

    \centering
    \caption{Execution time taken by the present algorithm relative to that taken by Algorithm~644 for computing the dataset described in Fig.~2 using double-precision arithmetic. Computations were performed using the five Fortran compilers described above on an Intel\textsuperscript{\textregistered} Core\textsuperscript{TM} i7-6600U CPU @ 2.60 GHz, 2.81 GHz processor.}
    \renewcommand{\arraystretch}{1.2}
    \begin{tabular}{lccccc}
        \hline
        \textbf{Region} & \(\boldsymbol{|z|\le 2.2}\)  &\(\boldsymbol{2.2<|z|\le 625}\) & \textbf{As. $\nu \to \infty$} &\(\boldsymbol{625<|z|}\) & \textbf{All} \\
        \hline
        \multicolumn{6}{l}{\textbf{GNU Fortran (i686-posix-dwarf-rev0, Built by MinGW-W64 project) 8.1.0}} \\
        Gfortran -O3 & 46.9\% & 53.7\% & 40.5\% & 40.5\% & 49.0\% \\
        Gfortran -O2 & 47.1\% & 52.2\% & 39.2\% & 39.2\% & 48.9\% \\
        Gfortran -O1 & 47.3\% & 54.6\% & 39.0\% & 39.0\% & 48.9\% \\
        Gfortran -O0 & 47.0\% & 55.4\% & 36.9\% & 36.9\% & 49.5\% \\
        \hline
        \multicolumn{6}{l}{\textbf{GNU Fortran (Rev3, Built by MSYS2 Project) 12.1.0}} \\
        Gfortran -O3 & 49.3\% & 91.4\% & 60.3\% & 70.3\% & 69.8\% \\
        Gfortran -O2 & 49.9\% & 89.2\% & 60.7\% & 67.6\% & 68.9\% \\
        Gfortran -O1 & 48.9\% & 85.0\% & 56.9\% & 75.7\% & 69.2\% \\
        Gfortran -O0 & 49.9\% & 75.2\% & 47.1\% & 73.1\% & 67.4\% \\
        \hline
        \multicolumn{6}{l}{\textbf{NAG Fortran Compiler Release 7.1 (Hanzomon) Build 7110}} \\
        Nagfor -O3 & 55.6\% & 80.8\% & 39.8\% & 75.1\% & 59.8\% \\
        Nagfor -O2 & 55.4\% & 96.0\% & 39.1\% & 96.7\% & 58.5\% \\
        Nagfor -O1 & 53.8\% & 56.6\% & 37.4\% & 95.1\% & 57.6\% \\
        Nagfor -O0 & 51.5\% & 58.7\% & 31.6\% & 56.7\% & 54.8\% \\
        \hline
        \multicolumn{6}{l}{\textbf{Intel(R) Fortran Intel(R) 64 Version 2021.9.0 Build 20230302\_000000}} \\
        Ifort -O3 & 62.2\% & 71.7\% & 47.1\% & 71.5\% & 63.8\% \\
        Ifort -O2 & 62.4\% & 69.1\% & 42.2\% & 67.2\% & 63.9\% \\
        Ifort -O1 & 62.4\% & 64.3\% & 42.4\% & 56.1\% & 63.3\% \\
        Ifort -O0 & 63.3\% & 69.7\% & 42.2\% & 71.5\% & 63.8\% \\
        \hline
        \multicolumn{6}{l}{\textbf{IFX (LLVM-based), from Intel® Fortran 64 Version 2021.9.0}} \\
        IFX -O3 & 63.7\% & 63.5\% & 45.5\% & 63.4\% & 61.2\% \\
        IFX -O2 & 63.4\% & 70.1\% & 43.5\% & 59.4\% & 61.4\% \\
        IFX -O1 & 64.1\% & 62.6\% & 43.5\% & 67.2\% & 61.0\% \\
        IFX -O0 & 63.7\% & 62.6\% & 47.2\% & 60.0\% & 60.2\% \\
        \hline
    \end{tabular}
    \label{tab:execution_time_comparison}
\end{table}

While the present algorithm exhibits superior performance across most of the computational domain, both in overall execution time and within specific subregions, there are small regions or isolated points where Algorithm~644 runs faster. Some of these occur on or near the overflow boundary (extremely large \(\nu\) and \(|z|\)), where Algorithm~644 may effectively bypass computations. Other instances correspond to high-order, high-argument values and are only marginally faster than the present algorithm, with no meaningful impact on overall performance.

As observed in Fig  \ref{fig:pt_by_pt_b}, the efficiency improvement is particularly pronounced in the small to intermediate orders where the present implementation adopts several efficiency improvements tacts such as those explained in applying Miller's algorithim in section \ref {sec:intermediate}. In addition, a noticeable improvement in efficiency can also be observed in the small $z$ region and the \(| z | \to \infty \) region.

Although the present algorithm and Algorithm~644 rely mainly on related analytical 
representations, such as power-series expansions, asymptotic expansions, and 
recurrence relations, the observed reduction in computational time arises from 
both differences in algorithmic design and implementation-level optimizations.

In particular, the present method reduces the region over which forward 
recurrence and Miller-type intermediate-region calculations are required. 
Large portions of the domain are instead treated using the asymptotic expansion 
for $z \to \infty$ and the uniform asymptotic expansion for $\nu \to \infty$, 
both of which are computationally more efficient when applicable.

Within the intermediate region, further efficiency improvements are introduced. 
For example, the term $\ln(C)$ in the estimate of the truncation parameter is 
replaced by an empirically chosen constant, and for sufficiently large $|z|$, 
fixed values of $M$ are used instead of repeatedly evaluating the full expression. 

In addition, the implementation avoids redundant evaluations of expensive 
functions. In particular, exponential terms are computed once and reused, 
rather than repeatedly evaluating equivalent hyperbolic functions of complex 
arguments. This reduces computational overhead and improves performance.

These factors collectively reduce the number of floating-point operations and 
the associated overhead, thereby explaining the observed speedup relative to 
Algorithm~644.

It is worth noting that the average execution time per evaluation using quad precision is a few \(\mu\)s, which is approximately an order of magnitude slower than double precision.

\clearpage 

\section{Conclusions}\label{sec:conclusion}

This work presents a robust and highly efficient algorithm and \texttt{Fortran} implementation for the computation of the modified Bessel function of the second kind, \( K_{\nu}(z) \), supporting both double- and quad-precision arithmetic. In line with Part~I of this study, which addressed \( I_{\nu}(z) \), the present implementation overcomes the critical limitations of existing algorithms, by extending the computational domain into regions where other algorithms fail.  

Extensive validation against high-precision \texttt{Maple} computations confirms the accuracy and reliability of the present algorithm, with relative errors  $< 10^{-26}$ for the real part and  $ <10^{-28}$ for the imaginary part in quadruple-precision arithmetic. These results are consistent with the performance established in Part~I for $I_{\nu}(z)$, thereby underscoring the robustness of the overall computational framework.

Efficiency benchmarks further highlight the advantages of the new implementation. Across all tested compilers and optimization levels, execution times are consistently reduced compared to Algorithm~644, making the algorithm particularly suitable for large-scale simulations and high-performance computing applications. The inclusion of quad-precision arithmetic enhances numerical stability, a feature of increasing importance in disciplines such as plasma physics, astrophysics, nuclear engineering, and quantum mechanics.  
 
Taken together, Parts~I and II provide a comprehensive, high-performance framework for computing the modified Bessel functions $I_{\nu}(z)$ and $K_{\nu}(z)$. The implementation is well documented and modular, allowing straightforward adaptation to precision formats beyond double and quadruple, including single, half, octuple, double-double, and quad-double precision. Having developed the algorithms and \texttt{Fortran} implementations for $I_{\nu}(z)$ and $K_{\nu}(z)$, extending the approach to design algorithms and computer codes for the regular Bessel functions becomes a relatively easy task.

\section*{Acknowledgments}
This work is supported by UAE University UPAR research grants \textbf{G00004187 (2022)} and \textbf{G00004995 (2024)}. 
The first author gratefully acknowledges the generous hospitality and support provided by the S. G. Johnson and the Department of Mathematics at MIT while serving as a Visiting Professor during Spring 2025.
\\
\\
\clearpage  
\setcounter{algorithm}{0}
\renewcommand{\thealgorithm}{\arabic{algorithm}}
\section*{Appendix:}

%

\textbf{Pseudocodes}
\begin{algorithm}[H]
\footnotesize
\setlength{\algorithmicindent}{0.8em}
\caption{\textsc{K\_nu\_of\_z}: Top-level dispatcher for $K_\nu(z)$}
\label{alg:KNU}
\begin{algorithmic}[1]
\STATE \textbf{Input:} $\nu \in \mathbb{R}$, $z_{\mathrm{in}}\in\mathbb{C}$ with $-\pi < \arg z_{\mathrm{in}} \le \pi$
\STATE \textbf{Output:} $K_\nu(z_{\mathrm{in}})$ in \texttt{cbk}, status \texttt{ierr}

\STATE \texttt{ierr} $\gets 0$
\IF{$z_{\mathrm{in}} = 0$}
  \STATE \texttt{ierr} $\gets 99$, \texttt{cbk} $\gets \mathrm{NaN}$, \textbf{return}
\ENDIF

\STATE $\nu_a \gets |\nu|$ \hfill\% evenness: $K_{-\nu}(z)=K_\nu(z)$
\STATE $x \gets \Re(z_{\mathrm{in}})$, $y \gets \Im(z_{\mathrm{in}})$, $r \gets \sqrt{x^2 + y^2}$
\STATE $z \gets z_{\mathrm{in}}$, $\mathrm{phase} \gets 1$, $s_y \gets 1$, \texttt{is\_left} $\gets \texttt{false}$

\IF{$x < 0$}
  \STATE \texttt{is\_left} $\gets \texttt{true}$, $z \gets -z_{\mathrm{in}}$
  \IF{$y = 0$}
    \STATE $s_y \gets 1$
  \ELSE
    \STATE $s_y \gets y/|y|$
  \ENDIF
  \STATE $n \gets \lfloor \nu_a \rfloor$, $f \gets \nu_a - n$, $\alpha \gets f\pi$
  \STATE $\mathrm{phase} \gets (-1)^n\bigl(\cos\alpha - i\,s_y\sin\alpha\bigr)$
\ENDIF

\STATE \texttt{series\_border} $\gets 2.2$, \texttt{z\_inf\_border} $\gets 625.0$, \texttt{inf\_nu\_const} $\gets 25.0$

\IF{$r \le \texttt{series\_border}$}
  \STATE call \textsc{K\_nu\_small\_z\_series}$(\nu_a, z)$ to obtain $(\texttt{cbk}, \texttt{ierr})$
\ELSIF{$\nu_a > \texttt{inf\_nu\_const} + r$}
  \STATE call \textsc{K\_nu\_unif\_sum}$(\nu_a, z)$ to obtain $(\texttt{cbk}, \texttt{ierr})$
\ELSIF{$r \le \texttt{z\_inf\_border}$}
  \STATE call \textsc{K\_nu\_intrmed\_z}$(\nu_a, z)$ to obtain $(\texttt{cbk}, \texttt{ierr})$
\ELSE
  \IF{$2r \ge \nu_a^2$}
    \STATE call \textsc{K\_nu\_inf\_z}$(\nu_a, z)$ to obtain $(\texttt{cbk}, \texttt{ierr})$
  \ELSE
    \STATE $n_{\mathrm{tgt}} \gets \lfloor \nu_a \rfloor$, $f \gets \nu_a - n_{\mathrm{tgt}}$
    \STATE $n_{\mathrm{seed}} \gets \lfloor \sqrt{2r} \rfloor - 2$, $n_{\mathrm{seed-1}} \gets n_{\mathrm{seed}} - 1$
    \STATE $\nu_{\mathrm{seed}} \gets n_{\mathrm{seed}} + f$, $\nu_{\mathrm{seed-1}} \gets n_{\mathrm{seed-1}} + f$
    \STATE call \textsc{K\_nu\_inf\_z}$(\nu_{\mathrm{seed}}, z)$ to obtain $(\texttt{Kseed}, \texttt{ierr1})$
    \STATE call \textsc{K\_nu\_inf\_z}$(\nu_{\mathrm{seed-1}}, z)$ to obtain $(\texttt{Kseedm1}, \texttt{ierr2})$
    \STATE \texttt{ierr} $\gets \max(\texttt{ierr1},\texttt{ierr2})$, $\texttt{two\_over\_z} \gets 2/z$
    \STATE $K_{n-1} \gets \texttt{Kseedm1}$, $K_n \gets \texttt{Kseed}$
    \FOR{$k = n_{\mathrm{seed}}$ \textbf{to} $n_{\mathrm{tgt}} - 1$}
      \STATE $v_{\mathrm{eff}} \gets f + k$
      \STATE $K_{n+1} \gets \texttt{two\_over\_z}\,v_{\mathrm{eff}} K_n + K_{n-1}$
      \STATE $K_{n-1} \gets K_n$, $K_n \gets K_{n+1}$
    \ENDFOR
    \STATE \texttt{cbk} $\gets K_n$
  \ENDIF
\ENDIF

\IF{\texttt{is\_left}}
  \STATE call \textsc{I\_nu\_of\_z}$(\nu_a, z)$ to obtain $(\texttt{cbi}, \texttt{ierrI})$
  \STATE \texttt{ierr} $\gets \max(\texttt{ierr},\texttt{ierrI})$
  \STATE \texttt{cbk} $\gets \mathrm{phase}\cdot\texttt{cbk} - i\,s_y\pi\,\texttt{cbi}$
\ENDIF

\STATE \textbf{return}
\end{algorithmic}
\end{algorithm}

\begin{algorithm}[H]
\footnotesize
\setlength{\algorithmicindent}{0.8em}
\caption{\textsc{K\_nu\_inf\_z}: Large-$|z|$ asymptotic for $K_\nu(z)$}
\label{alg:KINFZ}
\begin{algorithmic}[1]
\STATE \textbf{Input:} $\nu \ge 0$, $z\in\mathbb{C}$
\STATE \textbf{Output:} $K_\nu(z)$ in \texttt{cbk}, status \texttt{ierr}

\STATE \texttt{ierr} $\gets 0$
\STATE $x \gets \Re(z)$, $y \gets \Im(z)$, $r \gets \sqrt{x^2 + y^2}$
\STATE $\texttt{re\_lgK} \gets -x + \log\!\bigl(\sqrt{\pi/2}\bigr) - \tfrac12 \log r$
\IF{$\texttt{re\_lgK} < \texttt{log\_rmin}$}
  \STATE \texttt{ierr} $\gets -1$, \texttt{cbk} $\gets 0$, \textbf{return}
\ENDIF

\STATE $\texttt{two\_nu} \gets 2\nu$
\IF{$\texttt{two\_nu} > \texttt{sqrt\_rmin}$}
  \STATE $\texttt{four\_nu2} \gets \texttt{two\_nu}^2$
\ELSE
  \STATE $\texttt{four\_nu2} \gets 0$
\ENDIF

\STATE $\texttt{eight\_z} \gets 8z$, $\texttt{abs\_eight\_z} \gets 8r$
\STATE $\texttt{num} \gets 4\nu^2 - 1$
\STATE $\texttt{den} \gets \texttt{eight\_z}$, $\texttt{abs\_den} \gets \texttt{abs\_eight\_z}$
\STATE $\texttt{atolSq} \gets \left(\texttt{eps}\,\lvert\texttt{num}\rvert / \texttt{abs\_den}\right)^2$

\IF{$y = 0$ \textbf{and} $x > 0$}
  \STATE $\texttt{pref} \gets \sqrt{\pi/2}/\sqrt{x}$
  \STATE $\texttt{sum} \gets 1$, $\texttt{term} \gets 1$
  \FOR{$j = 1$ \textbf{to} $j_{\max}$}
    \STATE $\texttt{term} \gets \texttt{term}\cdot\texttt{num}/\texttt{den}$
    \STATE $\texttt{sum} \gets \texttt{sum}+\texttt{term}$
    \STATE $\texttt{den} \gets \texttt{den} + \texttt{eight\_z}$
    \STATE $\texttt{abs\_den} \gets \texttt{abs\_den} + \texttt{abs\_eight\_z}$
    \STATE $\texttt{num} \gets \texttt{num} - 8j$
    \IF{$\lvert\texttt{term}\rvert^2 \le \texttt{atolSq}$}
      \STATE \textbf{break}
    \ENDIF
  \ENDFOR
  \STATE \texttt{cbk} $\gets \texttt{pref}\,e^{-x}\,\texttt{sum}$
  \STATE \textbf{return}
\ENDIF

\STATE $\texttt{expFactor} \gets e^{-x}\bigl(\cos(-y) + i\sin(-y)\bigr)$
\STATE $s \gets \sqrt{z}$
\STATE $\texttt{pref} \gets \sqrt{\pi/2}/s$
\STATE $\texttt{sum} \gets 1$, $\texttt{term} \gets 1$

\FOR{$j = 1$ \textbf{to} $j_{\max}$}
  \STATE $\texttt{term} \gets \texttt{term}\cdot\texttt{num}/\texttt{den}$
  \STATE $\texttt{sum} \gets \texttt{sum}+\texttt{term}$
  \STATE $\texttt{den} \gets \texttt{den} + \texttt{eight\_z}$
  \STATE $\texttt{abs\_den} \gets \texttt{abs\_den} + \texttt{abs\_eight\_z}$
  \STATE $\texttt{num} \gets \texttt{num} - 8j$
  \IF{$\lvert\texttt{term}\rvert^2 \le \texttt{atolSq}$}
    \STATE \textbf{break}
  \ENDIF
\ENDFOR

\STATE \texttt{cbk} $\gets \texttt{pref}\cdot\texttt{expFactor}\cdot\texttt{sum}$
\STATE \textbf{return}
\end{algorithmic}
\end{algorithm}

\begin{algorithm}[H]
\footnotesize
\setlength{\algorithmicindent}{0.8em}
\caption{\textsc{K\_nu\_small\_z\_series}: Small-$|z|$ series for $K_\nu(z)$}
\label{alg:KSMALLZ}
\begin{algorithmic}[1]
\STATE \textbf{Input:} $\nu \ge 0$, $z\in\mathbb{C}$ with $|z| \le \texttt{series\_border}$
\STATE \textbf{Output:} $K_\nu(z)$ in \texttt{cbk}, status \texttt{ierr}

\STATE \texttt{ierr} $\gets 0$
\STATE $x\gets\Re(z)$, $y\gets\Im(z)$, $|z|\gets\sqrt{x^2+y^2}$
\STATE $z^{-1}\gets 1/z$, $\nu_{\rm eff}\gets\nu$

\IF{$\nu > \varepsilon$}
  \STATE $\lg_K \gets \log(\sqrt{\pi/2}) + (\nu-\tfrac12)\log\nu - \nu\bigl(0.3068528 + \log|z|\bigr)$
  \IF{$\lg_K > \log(R_{\max})$}
    \STATE \texttt{ierr}\,$\gets 1$
    \STATE \textbf{return}
  \ENDIF
\ENDIF

\IF{$\nu \le \tfrac12$}
  \IF{$|\nu-\tfrac12|<\varepsilon$}
    \STATE \texttt{cbk} $\gets
        \begin{cases}
          e^{-x}(\cos y - i\sin y)\sqrt{\tfrac{\pi}{2}}\,\sqrt{z^{-1}}, & x<|\log(R_{\min})|,\\
          0, & \text{otherwise}
        \end{cases}$
    \STATE \textbf{return}
  \ENDIF
  \STATE $\alpha\gets\nu$
\ELSE
  \STATE $n\gets\lfloor\nu\rfloor$, $\delta\gets\nu-n$
  \IF{$\delta=0$}
    \STATE $\alpha\gets0$, $n\gets n-1$
  \ELSIF{$\delta\ge\tfrac12$}
    \STATE $\alpha\gets\delta-1$
  \ELSE
    \STATE $\alpha\gets\delta$, $n\gets n-1$
  \ENDIF
\ENDIF

\STATE $(z/2)\gets \tfrac12 z$, $(z/2)^2\gets(z/2)^2$, $\alpha^2\gets\alpha^2$

\IF{$\alpha=0$}
  \STATE use limiting values for $\nu\to0$
\ELSE
  \STATE compute $g_1,g_2,q_0,p_0$ from $\Gamma(1\pm\alpha)$ and $(z/2)^{\pm\alpha}$
\ENDIF

\IF{$\nu < \tfrac12$}
  \STATE $c_0 \gets 1$, $f_0 \gets f_0(g_1,g_2,q_0,p_0)$, $\texttt{cbk} \gets c_0 f_0$
  \FOR{$k = 1$ \textbf{to} $k_{\max}$}
    \STATE update $f_k$, $c_k \gets c_{k-1}(z/2)^2/k$, $\Delta \gets c_k f_k$
    \STATE $\texttt{cbk} \gets \texttt{cbk} + \Delta$
    \IF{$|\Delta|^2 \le \varepsilon^2$}
      \STATE \textbf{break}
    \ENDIF
  \ENDFOR

\ELSIF{$\nu \in (\tfrac12,\tfrac32]$}
  \IF{$|\nu - \tfrac32| < \varepsilon$}
    \STATE \texttt{cbk} $\gets
      \begin{cases}
        e^{-x}(\cos y - i\sin y)\sqrt{\tfrac{\pi}{2}}\,\sqrt{z^{-1}}(1+z^{-1}), & x < |\log(R_{\min})|,\\
        0, & \text{otherwise}
      \end{cases}$
    \STATE \textbf{return}
  \ENDIF
  \STATE construct $K_\alpha$ by the above series; use $\nu$-recurrence to reach $K_\nu$

\ELSE
  \STATE compute $K_\alpha$ and $K_{\alpha+1}$ via the same series
  \FOR{$\mu=\alpha$ \textbf{to} $\nu$}
    \STATE $K_{\mu+1}(z) \gets K_{\mu-1}(z) + \bigl(2\mu/z\bigr)K_\mu(z)$
  \ENDFOR
  \STATE \texttt{cbk} $\gets K_\nu(z)$
\ENDIF

\STATE \textbf{return}
\end{algorithmic}
\end{algorithm}

\begin{algorithm}[H]
\footnotesize
\setlength{\algorithmicindent}{0.8em}
\caption{\textsc{K\_nu\_intrmed\_z}: Intermediate-$|z|$ region for $K_\nu(z)$}
\label{alg:KINTER}
\begin{algorithmic}[1]
\STATE \textbf{Input:} $\nu \ge 0$, $z\in\mathbb{C}$ with $\texttt{series\_border} < |z| \le \texttt{z\_inf\_border}$
\STATE \textbf{Output:} $K_\nu(z)$ in \texttt{cbk}, status \texttt{ierr}

\STATE \texttt{ierr} $\gets 0$, $x \gets \Re(z)$, $y \gets \Im(z)$, $r \gets \sqrt{x^2 + y^2}$

\IF{$\nu > \tfrac12$}
  \STATE $\ln K \gets \log(\sqrt{\pi/2}) + (\nu-\tfrac12)\log\nu - \nu\bigl(0.3068528+\log r\bigr)$
  \IF{$\ln K > \log(R_{\max})$}
    \STATE \texttt{ierr} $\gets 1$
    \STATE \textbf{return}
  \ENDIF
\ENDIF

\STATE $z^{-1} \gets 1/z$, $2z \gets 2z$

\IF{$r < \texttt{border1}$}
  \STATE $\theta \gets \arctan2(|y|,x)$, $a \gets 3\theta/(1+r)$, $b \gets 14.7\,\theta/(28+r)$
  \STATE $\log C \gets 32 + 24\cdot\texttt{rk\_by\_qp}$
  \STATE $\texttt{bracket} \gets [\log C + r\cos a/(1+0.008r)]/\cos b$
  \STATE $M \gets \Bigl(\tfrac{0.97}{2r}\,\texttt{bracket}^2 + 1.5\Bigr)\cdot\texttt{factor1}$
\ELSE
  \STATE $M \gets 7$ for double precision, and $M \gets 16$ for quad precision
\ENDIF

\IF{$\nu \le \tfrac12$}
  \STATE $\alpha \gets \nu$, $n \gets 0$
\ELSE
  \STATE $n \gets \lfloor \nu \rfloor$, $\delta \gets \nu - n$
  \IF{$\delta = 0$}
    \STATE $\alpha \gets 0$, $n \gets n-1$
  \ELSIF{$\delta \ge \tfrac12$}
    \STATE $\alpha \gets \delta - 1$
  \ELSE
    \STATE $\alpha \gets \delta$, $n \gets n-1$
  \ENDIF
\ENDIF

\STATE $q \gets \tfrac14 - \alpha^2$, $K_{M+1} \gets 0$, $K_M \gets \varepsilon$, $\texttt{sum} \gets K_{M+1}+K_M$

\FOR{$m = M$ \textbf{down to} $1$}
  \STATE $a_m \gets m^2 - m + q$, $b_m \gets \dfrac{2m + 2z}{m+1}$
  \STATE $K_{m-1} \gets \dfrac{-K_{m+1} + b_m K_m}{a_m}\,(m^2 + m)$
  \STATE $\texttt{sum} \gets \texttt{sum} + K_{m-1}$, $K_{m+1} \gets K_m$, $K_m \gets K_{m-1}$
\ENDFOR

\STATE $P(z) \gets \sqrt{\tfrac{\pi}{2z}}\,e^{-z}$
\IF{$x < |\log(R_{\min})|$}
  \STATE \texttt{cbk} $\gets P(z)\,K_0/\texttt{sum}$
\ELSE
  \STATE \texttt{cbk} $\gets 0$
\ENDIF

\IF{$n > 0$}
  \STATE compute $K_{\alpha+1}(z)$ from $K_\alpha(z)$ and the last recurrence
  \STATE $2/z \gets 2z^{-1}$, $K_\alpha \gets \texttt{cbk}$, $K_{\alpha+1}$ as above
  \FOR{$k = 1$ \textbf{to} $n$}
    \STATE $\mu \gets \alpha + k$, $K_{\mu+1} \gets (2\mu/z)K_\mu + K_{\mu-1}$
    \STATE update $K_{\mu-1}$ and $K_\mu$
  \ENDFOR
  \STATE \texttt{cbk} $\gets K_\nu$
\ENDIF

\STATE \textbf{return}
\end{algorithmic}
\end{algorithm}

\begin{algorithm}[H]
\footnotesize
\setlength{\algorithmicindent}{0.8em}
\caption{\textsc{K\_nu\_unif\_sum\_core}: Large-$\nu$ uniform asymptotic for $K_\nu(z)$}
\label{alg:KUNIF}
\begin{algorithmic}[1]
\STATE \textbf{Input:} large $\nu > 0$, $z\in\mathbb{C}$
\STATE \textbf{Output:} $K_\nu(z)$ in \texttt{cbk}, status \texttt{ierr}

\STATE This routine implements the uniform Debye--Olver expansion for large $\nu$ described in Part~I.
\STATE It uses the scaled variable $t = z/\nu$, auxiliary functions $p(t)$ and $\eta(t)$, and precomputed Debye coefficients in a convergent series in $1/\nu$.
\STATE Algorithm identical to the \textsc{K\_nu\_unif\_sum} routine in Part~I.
\STATE \textbf{return}
\end{algorithmic}
\end{algorithm}

\begin{algorithm}[H]
\footnotesize
\setlength{\algorithmicindent}{0.8em}
\caption{\textsc{I\_nu\_of\_z}: Dispatcher for $I_\nu(z)$ (Part~I)}
\label{alg:INUPARTI}
\begin{algorithmic}[1]
\STATE \textbf{Input:} $\nu\in\mathbb{R}$, $z\in\mathbb{C}$
\STATE \textbf{Output:} $I_\nu(z)$ in \texttt{cbi}, status \texttt{ierr}

\STATE This is the main $I_\nu(z)$ dispatcher developed in Part~I.
\STATE It selects among power series, large-$|z|$ asymptotics, uniform large-$\nu$ expansions, and Miller backward recurrence in $\nu$.
\STATE Implementation given in Part~I.
\STATE \textbf{return}
\end{algorithmic}
\end{algorithm}

\bibliographystyle{siamplain}
\bibliography{references}

\end{document}